\documentclass[reqno]{amsart}
\usepackage{amsthm,amsmath,amssymb}
\usepackage{stmaryrd,mathrsfs}

\newcommand{\ud}[0]{\,\mathrm{d}}

\newcommand{\abs}[1]{|#1|}
\newcommand{\Babs}[1]{\Big|#1\Big|}

\newcommand{\Norm}[2]{\|#1\|_{#2}}

\newcommand{\BNorm}[2]{\Big\|#1\Big\|_{#2}}
\newcommand{\pair}[2]{\langle #1,#2 \rangle}
\newcommand{\Bpair}[2]{\Big\langle #1,#2 \Big\rangle}
\newcommand{\ave}[1]{\langle #1\rangle}
\newcommand{\diag}[0]{\operatorname{diag}}

\newcommand{\bddlin}[0]{\mathscr{L}}

\newcommand{\R}{\mathbb{R}}
\newcommand{\C}{\mathbb{C}}

\newcommand{\Sb}{\mathbb{S}}

\newcommand{\Exp}[0]{\mathbb{E}}

\newcommand{\ontop}[2]{\begin{smallmatrix} #1 \\ #2 \end{smallmatrix}}

\numberwithin{equation}{section}

\makeatletter
  \let\c@equation\c@subsection
\makeatother

\theoremstyle{plain}

\newtheorem{theorem}[subsection]{Theorem}
\newtheorem{proposition}[subsection]{Proposition}

\theoremstyle{definition}

\theoremstyle{remark}

\begin{document}

\title[Beurling--Ahlfors operator]{On the norm of the Beurling--Ahlfors operator in several dimensions}

\author[T.~P.\ Hyt\"onen]{Tuomas P.\ Hyt\"onen}
\address{Department of Mathematics and Statistics, University of Helsinki, Gustaf H\"allstr\"omin katu 2b, FI-00014 Helsinki, Finland}
\email{tuomas.hytonen@helsinki.fi}



\maketitle

\begin{abstract}
The generalized Beurling--Ahlfors operator \(S\) on \(L^p(\R^n;\Lambda)\), where \(\Lambda:=\Lambda(\R^n)\) is the exterior algebra with its natural Hilbert space norm, satisfies the estimate
\begin{equation*}
  \|S\|_{\mathscr{L}(L^p(\R^n;\Lambda))}\leq(n/2+1)(p^*-1),\qquad p^*:=\max\{p,p'\}.
\end{equation*}
This improves on earlier results in all dimensions \(n\geq 3\). The proof is based on the heat extension and relies at the bottom on Burkholder's sharp inequality for martingale transforms.
\end{abstract}

\section{Introduction}

This note is an outgrowth of my study of Ba\~nuelos \& Lindeman's paper~\cite{BaLin}, which deals with the same theme: the estimation of the norm of the generalized Beurling--Ahlfors operator \(S\) which acts on the exterior algebra \(\Lambda:=\Lambda(\R^n)\) -valued functions in arbitrary dimension \(n\geq 2\). This operator, or more precisely the generalization of the classical planar version to \(n>2\), was introduced by Iwaniec \& Martin~\cite{IwaM}, who also established the bounds
\begin{equation*}
  (p^*-1)\leq\Norm{S}{\bddlin(L^p(\R^n;\Lambda))}
  \leq C(n)(p^*-1),\qquad p^*=\max\{p,p'\},
\end{equation*}
and conjectured that in fact \(C(n)=1\) in all dimensions \(n\geq 2\) (\cite{IwaM}, p.~34). For \(n=2\), both the theorem and the conjecture date further back, the lower bound being due to Lehto~\cite{Lehto}.

The problem remains open even in the plane, where the best upper bound at the time of writing appears to be that of Ba\~nuelos \& Janakiraman~\cite{BaJan} with \(C(2)\leq 1.575\). A wealth of motivation for such norm estimates can be found in the papers cited above and below in this introduction, so it seems redundant to repeat that discussion here.

In higher dimensions, the sharpest estimates until recently were those of Ba\~nuelos \& Lindeman~\cite{BaLin}, who obtained \(C(n)\leq 4n/3+O(1)\). These bounds rely on a representation of \(S\) as a transform of Brownian martingales arising from the harmonic extension of the function \(f\) on which \(S\) acts, and then on Burkholder's sharp inequality for martingale transforms~\cite{Burkholder}.

Since Ba\~nuelos \& Lindeman's paper~\cite{BaLin}, Nazarov \& Volberg~\cite{NaVo} have discovered that it is more efficient to use the heat extension than the harmonic extension in the estimation of the Beurling--Ahlfors operators: this idea, implemented as a Bellman function argument, improved the record bound at the time, \(C(2)\leq 4\) by Ba\~nuelos \& Wang~\cite{BaWan}, down to \(C(2)\leq 2\). Very recently, the heat--Bellman strategy was also employed in higher dimensions by Petermichl, Slavin \& Wick~\cite{PSW}, who obtained \(C(n)\leq n\).

In view of the basic difference between the harmonic and heat extension methods, one would expect to cut off half of the constant of Ba\~nuelos \& Lindeman~\cite{BaLin} also in the higher dimensional case, i.e., to get \(C(n)\leq 2n/3+O(1)\). Indeed, this can be reached by carefully following their original proof but with the harmonic extension replaced by the heat extension throughout, and this will be done in the present paper. Actually, I even obtain the better estimate \(C(n)\leq n/2+O(1)\), thanks to an additional elementary trick which exploits the flexibility resulting from the non-uniqueness of the heat representation of the Beurling--Ahlfors operator. I also derive bounds for the restrictions of \(S\) to \(r\)-forms, which similarly improve on the corresponding results of Ba\~nuelos \& Lindeman~\cite{BaLin}.

The precise statement is the following:

\begin{theorem}\label{thm:main}
The Beurling--Ahlfors operator \(S\), and its restriction to \(r\)-forms, satisfy the following norm bounds in all dimensions \(n\geq 2\), for all \(r\in\{0,1,\ldots,n\}\), and for all \(p\in(1,\infty)\):
\begin{equation*}\begin{split}
  \Norm{S}{\bddlin(L^p(\R^n;\Lambda^r))}
  &\leq\Big(\frac{2r(n-r)}{n}+1\Big)(p^*-1)
  \leq(2r+1)(p^*-1), \\
  \Norm{S}{\bddlin(L^p(\R^n;\Lambda))}
  &\leq(p^*-1)\times\begin{cases}
    (n/2+1),&n\text{ even}, \\
    (n/2+1-1/2n), &n\text{ odd}.
  \end{cases}
\end{split}\end{equation*}
\end{theorem}

The latter expression in the estimate on \(L^p(\R^n;\Lambda^r)\) shows that there is a dimen\-sion-free (i.e., \(n\)-free) bound for a fixed \(r\), as already discovered by Ba\~nuelos \& Lindeman~\cite{BaLin} with a bigger constant, but as the middle expression tells, a better estimate is available by taking into account the dimension \(n\). As mentioned in the beginning, the bound of the Theorem is not new for \(n=2\), and indeed one can do better with the more refined methods of Ba\~nuelos \& Janakiraman~\cite{BaJan} in this case. In all higher dimensions \(n\geq 3\), the result improves on previous estimates, as far as I am aware of.

In the planar case, there are also asymptotic estimates for the Beurling--Ahlfors operator as \(p\to\infty\), which come somewhat closer to the conjecture than the results known for arbitrary \(p\in(1,\infty)\). More precisely, Dragi\v{c}evi\'c \& Volberg~\cite{DraVo} showed that \(\Norm{S}{\bddlin(L^p(\R^2;\Lambda))}\leq\sqrt{2}(1+o(1))(p-1)\) as \(p\to\infty\). Since this result was based on an asymptotic refinement of Nazarov \& Volberg's~\cite{NaVo} \(C(2)\leq 2\) result, one might optimistically hope that a similar approach could provide \(n\)-dimensional large-\(p\) asymptotics of the order \(O(\sqrt{n})(p-1)\). In this light the following result that I actually managed to show is perhaps not too impressive. Nevertheless it seems worth proving, if only for the sake of revisiting Dragi\v{c}evi\'c \& Volberg's nice argument and perhaps clarifying the essence of the underlying idea by generalizing it to several variables.

\begin{proposition}\label{prop:asymp}
The Beurling--Ahlfors operator \(S\) satisfies the following asymptotic norm bounds in all dimensions \(n\geq 2\) as \(p\to\infty\):
\begin{equation*}
  \Norm{S}{\bddlin(L^p(\R^n;\Lambda))}
  \leq(1+o(1))(p-1)\times\begin{cases}
    \sqrt{(n/2)^2+1},&n\text{ even}, \\
    \sqrt{(n/2)^2+3/4}, &n\text{ odd}.
  \end{cases}
\end{equation*}
\end{proposition}

Theorem~\ref{thm:main} and Proposition~\ref{prop:asymp} are proved in the rest of the paper. In the final section I also provide some results of the same flavour for the spectral multipliers of the Laplace operator, its imaginary powers in particular.

\section{The heat matrix representing an operator}

Let \(L^p(\R^n;\R^N)\) stand for the Lebesgue space of \(\R^N\)-valued \(p\)-integrable functions \(f(x)=\{f_j(x)\}_{j=1}^N\) on \(\R^n\) with the usual norm
\begin{equation*}
  \Norm{f}{L^p(\R^n;\R^N)}:=
  \Big(\int_{\R^n}\Big[\sum_{j=1}^N\abs{f_j(x)}^2\Big]^{p/2}\ud x\Big)^{1/p}.
\end{equation*}
Suppose that an operator \(T\in\bddlin(L^p(\R^n;\R^N),L^p(\R^n;\R^M))\) has a representation
\begin{equation}\label{eq:TvsA}
  \int_{\R^n}\pair{Tf(x)}{g(x)}\ud x
  =\int_0^{\infty}\int_{\R^n}\pair{A(x,t)\nabla u(x,t)}{\nabla v(x,t)}\ud x\ud t,
\end{equation}
where \(\nabla\) refers to the \(x\)-gradient (the gradient of an \(\R^N\)-valued function being an \(\R^{nN}=\R^n\otimes\R^N\)-valued one), \(u(x,t)=e^{\frac{1}{2}t\triangle}f(x)\) and \(v(x,t)=e^{\frac{1}{2}t\triangle}g(x)\) are the \emph{heat extensions} of \(f\in L^p(\R^n;\R^N)\) and \(g\in L^p(\R^n;\R^M)\), and \(A(x,t)\in\bddlin(\R^{nN},\R^{nM})\) is an \(nN\times nM\)-matrix for all \((x,t)\in\R_+^{n+1}\). Then it is a by-now well-known fact (cf. Ba\~nuelos \& M\'endez-Hern\'andez~\cite{BaMen} for \(n=2\), \(N=M=1\)) that
\begin{equation}\label{eq:Burkholder}
  \Norm{T}{\bddlin(L^p(\R^n;\R^N),L^p(\R^n;\R^M))}
  \leq\sup_{(x,t)\in\R_+^{n+1}}\Norm{A(x,t)}{\bddlin(\R^{nN},\R^{nM})}(p^*-1).
\end{equation}

\begin{proof}
I give this for the convenience of the reader, although it differs little from the related arguments found in the literature. Assume without loss of generality that \(\Norm{A(x,t)}{\bddlin(\R^{nN},\R^{nM})}\leq 1\). It clearly suffices to show that
\begin{equation}\label{eq:BurkReduct}\begin{split}
  &\lim_{\tau\to\infty}\Babs{\int_0^{\tau}\int_{\R^n}\pair{A(x,\tau-t)\nabla u(x,\tau-t)}{\nabla v(x,\tau-t)}\ud x\ud t} \\
  &\qquad\leq(p^*-1)\Norm{f}{L^p(\R^n;\R^N)}\Norm{g}{L^{p'}(\R^n;\R^M)}
\end{split}\end{equation}
This is proved by a probabilistic argument.

Let \(\Exp^x\) stand for the expectation related to a probability measure governing an \(\R^n\)-valued standard Brownian motion \((X_t)_{t\in[0,\tau]}\) starting at \(x\in\R^n\). Then by the Markov property
\begin{equation}\label{eq:Markov}
  \int_{\R^n}\Exp^x g(X_t)\ud x
  =\int_{\R^n}e^{\frac{1}{2}t\triangle}g(x)\ud x
  =\int_{\R^n}g(x)\ud x,
\end{equation}
valid for any \(g\in L^1(\R^n)\), and the It\^o isometry for stochastic integrals,
\begin{equation}\label{eq:ItoIsom}\begin{split}
  &\int_0^{\tau}\int_{\R^n}\pair{A(x,\tau-t)\nabla u(x,\tau-t)}{\nabla v(x,\tau-t)}\ud x\ud t \\
  &=\int_0^{\tau}\int_{\R^n}\Exp^x\pair{A(X_t,\tau-t)\nabla u(X_t,\tau-t)}{\nabla v(X_t,\tau-t)}\ud x\ud t \\
  &=\int_{\R^n}\Exp^x\Bpair{\int_0^{\tau}A\nabla u(X_t,\tau-t)\cdot\ud X_t}{
    \int_0^{\tau}\nabla v(X_t,\tau-t)\cdot\ud X_t}\ud x,
\end{split}\end{equation}
where I have abbreviated \(A\nabla u(x,t):=A(x,t)\nabla u(x,t)\).

The right side of \eqref{eq:ItoIsom} requires some interpretation, as \(Au\) and \(v\) are \(nN\)- and \(nM\)-vectors, respectively, while \(X_t\) is an \(n\)-vector: it is understood that the dot products are taken only with respect to the coordinates in \(\R^n\), i.e., for \(\xi=\{\xi_{i,j}\}_{1\leq i\leq n;1\leq j\leq N}\in\R^{nN}\) and \(\eta=\{\eta_i\}_{1\leq i\leq n}\in\R^n\), the dot product is \(\xi\cdot\eta=\{\sum_{i=1}^n\xi_{i,j}\eta_i\}_{1\leq j\leq N}\in\R^N\). Next, by H\"older,
\begin{equation*}\begin{split}
  \abs{RHS\eqref{eq:ItoIsom}}
  \leq &\Big(\int_{\R^n}\Exp^x\BNorm{\int_0^{\tau}A\nabla u(X_t,\tau-t)\cdot\ud X_t}{\R^M}^p\ud x\Big)^{1/p} \\
     &\times\Big(\int_{\R^n}\Exp^x\BNorm{\int_0^{\tau}\nabla v(X_t,\tau-t)\cdot\ud X_t}{\R^M}^{p'}\ud x\Big)^{1/p'}.
\end{split}\end{equation*}

Let \(Y_s:=\int_0^{s}A\nabla u(X_t,\tau-t)\cdot\ud X_t\) and \(U_s:=\int_0^s\nabla u(X_t,\tau-t)\cdot X_t\). Then \((Y_s)_{s\in[0,\tau]}\) and \((U_s)_{s\in[0,\tau]}\) are \(\R^M\)- and \(\R^N\)-valued continuous-path martingales, respectively. Their quadratic variations are
\begin{equation*}
  \ave{Y}_s=\int_0^s\Norm{A\nabla u(X_t,\tau-t)}{\R^{nM}}^2\ud t,\qquad
  \ave{U}_s=\int_0^s\Norm{\nabla u(X_t,\tau-t)}{\R^{nN}}^2\ud t.
\end{equation*}
Then \(\ave{U}_s-\ave{Y}_s\) is a non-negative increasing process, as it starts at \(0\) and the increment from \(r\) to \(s>r\) is given by
\begin{equation*}
  \int_r^s\Big(\Norm{\nabla u(X_t,\tau-t)}{\R^{nN}}^2
     -\Norm{A\nabla u(X_t,\tau-t)}{\R^{nM}}^2\Big)\ud t\geq 0.
\end{equation*}

Hence \(Y_s\) and \(U_s\) satisfy the assumptions of Ba\~nuelos and M\'endez-Hern\'andez~\cite{BaMen}, Theorem~2 (which in turn is based on related results due to Burkholder~\cite{Burkholder}). The mentioned theorem guarantees that
\begin{equation*}
  \big(\Exp^x\Norm{Y_{\tau}}{\R^M}^p\big)^{1/p}
  \leq(p^*-1)\big(\Exp^x\Norm{U_{\tau}}{\R^N}^p\big)^{1/p}.
\end{equation*}

Next one observes that
\begin{equation*}\begin{split}
  U_{\tau}
  &=\int_0^{\tau}\nabla u(X_t,\tau-t)\cdot\ud X_t \\
  &=u(X_{\tau},0)-u(X_0,\tau)-\int_0^{\tau}(-\partial_t u+\frac{1}{2}\triangle)u(X_t,\tau-t)\ud t \\
  &=u(X_{\tau},0)-u(X_0,\tau)=f(X_{\tau})-e^{\frac{1}{2}\tau\triangle}f(X_0)
\end{split}\end{equation*}
by It\^o's formula and the fact that \(u\) is the heat extension of \(f\). By the Markov property \eqref{eq:Markov},
\begin{equation*}
  \Big(\int_{\R^n}\Exp^x\Norm{f(X_{\tau})}{\R^N}^p\ud x\Big)^{1/p}
  =\Norm{f}{L^p(\R^n;\R^N)},
\end{equation*}
whereas
\begin{equation*}
  \Big(\int_{\R^n}\Exp^x\Norm{e^{\frac{1}{2}\tau\triangle}f(X_0)}{\R^N}^p\ud x\Big)^{1/p}
  =\Norm{e^{\frac{1}{2}\tau\triangle}f}{L^p(\R^n;\R^N)}
  \to 0\qquad\text{as }\tau\to\infty.
\end{equation*}
Repeating the last few steps with \((v,g,p',M)\) in place of \((u,f,p,N)\) completes the proof of \eqref{eq:Burkholder}.
\end{proof}

It might be interesting to remark that with the choice
\begin{equation*}
  A(x,t)=\frac{\nabla u(x,t)\otimes\nabla v(x,t)}{\Norm{\nabla u(x,t)}{\R^{nN}}\Norm{\nabla v(x,t)}{\R^{nM}}}
\end{equation*}
in \eqref{eq:TvsA}, the estimate \eqref{eq:Burkholder} gives
\begin{equation}\label{eq:PSW}\begin{split}
  &\int_0^{\infty}\int_{\R^n}
   \Norm{\nabla u(x,t)}{\R^{nN}}\Norm{\nabla v(x,t)}{\R^{nM}}\ud x\ud t \\
  &\qquad\leq(p^*-1)\Norm{f}{L^p(\R^n;\R^N)}\Norm{g}{L^{p'}(\R^n;\R^M)},
\end{split}\end{equation}
which is Petermichl, Slavin \& Wick's~\cite{PSW} Theorem~1.3, there proved by a Bellman function argument. (Their heat extension is defined with the semigroup \(e^{t\triangle}\) instead of \(e^{\frac{1}{2}t\triangle}\) used here, which explains the absence of the factor \(2\) on the left side of \eqref{eq:PSW} as compared to the statement in~\cite{PSW}.)

Conversely, \eqref{eq:Burkholder} is an immediate consequence of \eqref{eq:PSW}. This is not a coincidence, since despite their superficial dissimilarity, the Bellman function argument of~\cite{PSW} and the present proof are just reflections of the same underlying phenomenon, i.e., the fundamental martingale inequalities of Burkholder~\cite{Burkholder}.

\section{The case of the Beurling--Ahlfors operator}

Now consider the particular case of the previous section when \(N=M=2^n\) and \(\R^N\) is identified with the exterior algebra \(\Lambda\) with canonical unit vectors
\begin{equation*}
  e_I=e_{i_1}\wedge e_{i_2}\wedge\cdots\wedge e_{i_r},\qquad
  I=\{i_1,i_2,\ldots,i_r\},\quad 1\leq i_1<i_2<\ldots<i_r\leq n,
\end{equation*}
indexed by the subsets \(I\subseteq\{1,2,\ldots,n\}\). The operator of interest is the Beurling--Ahlfors transform \(S\), which will be here defined by the following formula established by Iwaniec \& Martin (\cite{IwaM}, top of p.~58):
\begin{equation}\label{eq:Sf}
  Sf=
  \sum_K\Big[\sum_{k\in K}R_k^2-\sum_{\ell\notin K}R_{\ell}^2\Big]f_K e_K
  +\sum_K\sum_{\ontop{k\in K}{\ell\notin K}}2R_k R_{\ell} f_K e_{K-k+\ell},
\end{equation}
where \(f=\sum_K f_K e_K\), the symbols \(R_1,\ldots,R_n\) stand for the usual Riesz transforms on \(\R^n\), and, given \(e_K=e_{i_1}\wedge\cdots\wedge e_k\wedge\cdots\wedge e_{i_r}\), the vector \(e_{K-k+\ell}\) is defined by \(e_{K-k+\ell}:=e_{i_1}\wedge\cdots\wedge e_{\ell}\wedge\cdots\wedge e_{i_r}\), i.e., by substituting \(e_{\ell}\) in place of \(e_k\). Note that \(e_{K-k+\ell}\) may fail to be one of the canonical unit vectors, since \(e_{\ell}\) may be in a ``wrong'' place, but with elementary algebra one checks that
\begin{equation}\label{eq:sign}
  e_{K-k+\ell}=(-1)^{\# K(k,\ell)}e_{K\setminus k\cup\ell},
\end{equation}
where \(K(k,\ell):=K\cap(\min\{k,\ell\},\max\{k,\ell\})\) is the part of \(K\) between \(k\) and \(\ell\). The indentification \(k=\{k\}\), \(\ell=\{\ell\}\) is made to simplify writing in the expression \(K\setminus k\cup\ell:=K\setminus\{k\}\cup\{\ell\}\). Apparently, the sign in \eqref{eq:sign} has been missed in some of the related literature.

A representation of the type \eqref{eq:TvsA} for \(S\) is readily derived by combining \(\eqref{eq:Sf}\) with the well-known representations of the Riesz transforms. Let \(E_{k\ell}\) designate the \(n\times n\)-matrix consisting of zeros except for a one at the position \((k,\ell)\). Similarly, \(E_{KL}\) (with capital indices representing sets) will stand for the \(2^n\times 2^n\)-matrix of zeros except for a one at the place \((K,L)\). Then \eqref{eq:TvsA} (with \(N=1\)) holds for \(T=-R_k R_{\ell}\) and \(A=\alpha E_{k\ell}+(1-\alpha)E_{\ell k}\), where the choice of \(\alpha\in\R\) is arbitrary. The freedom to pick any \(\alpha\in[0,1]\) will be exploited below. Because of the minus sign above, I will instead present the matrix \(A\) representing \(-S\), which is
\begin{equation}\label{eq:AofS}\begin{split}
  A &=\sum_K\Big[\sum_{k\in K}E_{kk}-\sum_{\ell\notin K}E_{\ell\ell}\Big]\otimes E_{KK} \\
  &+\sum_K\sum_{\ontop{k\in K}{\ell\notin K}}2[\alpha_{k\ell}(K) E_{k\ell}
  +(1-\alpha_{k\ell}(K)) E_{\ell k}](-1)^{\# K(k,\ell)}\otimes E_{K\setminus k\cup\ell,K}.
\end{split}\end{equation}
Note that \(A\) is just a constant matrix here, so the full generality of \eqref{eq:Burkholder} allowing dependence on \(x\) and \(t\) is not used in the application to the Beurling--Ahlfors operator.

It can be seen at once that \(A\) is of a block form, with interaction only between sets \(K\) of the same size. Fix a number \(r\in\{0,1,\ldots,n\}\), and consider the part of \(A\) with \(\# K=r\). As suggested in \eqref{eq:AofS}, there is in principle the freedom of choosing a different value of \(\alpha\) for every triplet \((k,\ell,K)\), but I will only exploit this by making \(\alpha\) a function of \(r\), so the indicated dependence on \((k,\ell,K)\) will be dropped in the subsequent analysis. There is a chance that a more sophisticated choice of the \(\alpha\)'s would lead to slightly better estimates, but I have decided not to take up this additional complication, since it seems that the possible improvement, if any, obtainable this way would be quite insignificant.

Let me mention that choosing \(\alpha\) as a function of \(r\) is the advertised elementary trick which allows to improve the estimates of Ba\~nuelos and Lindeman~\cite{BaLin} by more than the factor \(\frac{1}{2}\) which one would get by only repeating their argument with the heat extension in place of the harmonic extension. In fact, a restricted version of the same trick was already applied in~\cite{BaLin}, but only the three values \(\alpha\in\{0,\frac{1}{2},1\}\) were exploited, corresponding to the completely one-sided and completely symmetric representations of \(R_iR_j\).

\section{Block structure of the matrix}

I will next compute the entry of \(A\) in the position \((I,i;J,j)\), where \(i,j\) are elements and \(I,J\) subsets of \(\{1,2,\ldots,n\}\). Whenever \(\mathscr{P}\) is a logical statement, \([\mathscr{P}]\) will designate \(1\) if \(\mathscr{P}\) is true and \(0\) if \(\mathscr{P}\) is false. E.g., \([i\in I]\) is the same as the characteristic function \(1_I(i)\), but the bracket notation is sometimes more handy to avoid complicated trees of sub- and superscripts.
\begin{equation}\label{eq:Aij}\begin{split}
  A_{Ii,Jj} &=\sum_K\Big[\sum_{k\in K}\delta_{ik}\delta_{jk}
     -\sum_{\ell\notin K}\delta_{i\ell}\delta_{j\ell}\Big]\delta_{IK}\delta_{JK} \\
  &+\sum_K\sum_{\ontop{k\in K}{\ell\notin K}}2\big(\alpha\delta_{ik}\delta_{j\ell}
  +(1-\alpha)\delta_{i\ell}\delta_{jk}\big)(-1)^{\# K(k,\ell)}\delta_{I,K\setminus k\cup\ell}\delta_{JK} \\
  &=\delta_{IJ}\delta_{ij}(-1)^{[i\notin I]}
   +[i\in J,j\notin J,I=J\setminus i\cup j]\cdot 2\alpha(-1)^{\# J(i,j)} \\
  &\qquad+[i\notin J,j\in J,I=J\setminus j\cup i]\cdot 2(1-\alpha)(-1)^{\# J(i,j)}.
\end{split}\end{equation}

From \eqref{eq:Aij} one sees that if a non-diagonal entry \(A_{Ii,Jj}\) is non-zero, then either \(i\in I\) and \(j\in J\), or else \(i\notin I\) and \(j\notin J\), i.e., the pairs \((I,i)\) and \((J,j)\) are in this sense of the same type, and pairs of different type belong to different blocks.

Consider first the pairs of the type \(i\notin I\). Let \(\tilde{I}:=\{i_0,i_1,\ldots,i_r\}\) and \(I_0=\tilde{I}\setminus i_0\). One sees that \(A_{I_0,i_0;J,j}\) or \(A_{J,j;I_0,i_0}\) is non-zero only if \((J,j)=(I_k,i_k)\) for some \(k=0,1,\ldots,r\), where \(I_k:=\tilde{I}\setminus i_k\). Hence each \((I,i)\) of the type \(i\notin I\) belongs to a block of size \((r+1)\times(r+1)\), corresponding to the columns and rows of the form \((I_k,i_k)\) in the full matrix. It may be assumed that \(i_0<i_1<\ldots<i_r\), whence it is easy to see that, for \(s<t\), \(I_s(s,t)=I_t(s,t)=\{i_{s+1},\ldots,i_{t-1}\}\) has size \(t-s-1=\abs{t-s}-1\). Hence this block has \(-1\)'s on the diagonal and \(2\alpha(-1)^{s+t+1}\) in any position \((t,s)\) away from the diagonal, so it is equal to
\begin{equation}\label{eq:notin}\begin{split}
 &\begin{pmatrix}
     -1 & 2\alpha & -2\alpha & 2\alpha & \\
     2\alpha & -1 & 2\alpha & -2\alpha & \cdots \\
     -2\alpha & 2\alpha & -1  & 2\alpha & \\
    2\alpha & -2\alpha & 2\alpha & - 1 & \\
     & \vdots & & & \ddots \\
  \end{pmatrix} \\
 &=  2\alpha \begin{pmatrix}
     -1 & 1 & -1 & \\
     1 & -1 & 1 & \cdots \\
     -1 & 1 & -1 & \\
        & \vdots & & \ddots \\
   \end{pmatrix}
   +(2\alpha-1)I_{r+1},
\end{split}\end{equation}
where \(I_{r+1}\) is the \((r+1)\times(r+1)\) identity matrix. Let also \(J_{r+1}\) stand for the \((r+1)\times(r+1)\) matrix with all entries equal to one. It is clear that the matrix with alternating \(+1\)'s and \(-1\)'s appearing in the right of \eqref{eq:notin} is similar to \(-J_{r+1}\), the similarity being realized by the diagonal matrix \(\diag(1,-1,1,-1,\ldots)\). Hence the full matrix in \eqref{eq:notin} is similar to
\begin{equation*}\begin{split}
  -2\alpha J_{r+1}+(2\alpha-1)I_{r+1}
  &\sim -2\alpha\diag((r+1),0,\ldots,0)+(2\alpha-1)I_{r+1} \\
  &=\diag(-(2\alpha r+1),2\alpha-1,\ldots,2\alpha-1).
\end{split}\end{equation*}
Here \(\sim\) means similarity of matrices, and this is realized by the unitary matrix \(x_{ts}=\exp\big(i2\pi ts/(r+1)\big)/\sqrt{r+1}\), as can be easily checked. It follows that the norm of the block in \eqref{eq:notin} is
\begin{equation*}
  \max\{\abs{2\alpha r+1},\abs{2\alpha-1}\}=2\alpha r+1,
\end{equation*}
since \(r\geq 0\) and \(\alpha\in[0,1]\).

It remains to treat the pairs of the type \(i\in I\). Fix some \(\tilde{I}:=\{i_1,\ldots,i_{r-1}\}\) where \(i_1<i_2<\ldots<i_{r-1}\), and let \(I_i:=\tilde{I}\cup i\) for some \(i\notin\tilde{I}\). Then \(A_{I_i,i;J,j}\) or \(A_{J,j;I_i,i}\) is non-zero only if \(J\) is also of the form \(J=\tilde{I}\cup j=:I_j\) for some \(j\notin\tilde{I}\). Hence every pair \((I,i)\) of the type \(i\in I\) belongs to a block of size \(\#\tilde{I}^c=n-(r-1)=n-r+1\). The diagonal elements \(A_{Ii,Ii}\) are equal to \(1\). As for the signs of the non-diagonal elements, observe that
\begin{equation*}
  I_i(i,j)=I_j(i,j)=\tilde{I}(i,j)=\tilde{I}\cap(\min\{i,j\},\max\{i,j\}),
\end{equation*}
so the non-diagonal signs should be chosen according to the matrix
\begin{equation}\label{eq:signs}
  \begin{matrix}
             & 0<j<i_1 & i_1<j<i_2 & i_2 < j<i_3 & \cdots \\
   0<i<i_1     & + & - & + & \\
   i_1<i<i_2 & - & + & - & \cdots \\
   i_2<i<i_3 & + & - & + & \\
    \vdots   & & \vdots & & \ddots \\
  \end{matrix}
\end{equation}
Incidentally, this also gives the correct sign of the diagonal elements. Observe that some of the blocks \(i_s<j<i_{s+1}\) may be empty if \(\tilde{I}\) contains consecutive elements. Again, a matrix with signs as in \eqref{eq:signs} is similar to one with plus signs everywhere, the similarity being realized by a diagonal matrix of signs as in, say, the first line of \eqref{eq:signs}. So up to similarity the relevant block is now
\begin{equation*}\begin{split}
  &\begin{pmatrix}
     1 & 2(1-\alpha) & 2(1-\alpha) & \\
     2(1-\alpha) & 1 & 2(1-\alpha) & \cdots \\
     2(1-\alpha) & 2(1-\alpha) & 1 & \\
     & \vdots & & \ddots \\
  \end{pmatrix} \\
  &=2(1-\alpha)J_{n-r+1}+(2\alpha-1)I_{n-r+1} \\
  &\sim 2(1-\alpha)\diag(n-r+1,0,\ldots,0)+(2\alpha-1)I_{n-r+1} \\
  &=\diag(2(1-\alpha)(n-r)+1,2\alpha-1,\ldots,2\alpha-1).
\end{split}\end{equation*}
This has norm
\begin{equation*}
  \max\{\abs{2(1-\alpha)(n-r)+1},\abs{2\alpha-1}\}
  =2(1-\alpha)(n-r)+1,
\end{equation*}
since \(r\leq n\) and \(\alpha\in[0,1]\).

\section{Norm of the matrix}

The outcome of the previous section was to quantify the norm of the part of \(A\) related to \(r\)-sets as
\begin{equation}\label{eq:partNorm}
  \max\{2\alpha r+1,2(1-\alpha)(n-r)+1\},
\end{equation}
where the choice of \(\alpha\in[0,1]\) is still to be made. The quantity in \eqref{eq:partNorm} will be minimized by requiring the two expressions to be equal. This yields \(\alpha=1-r/n\), which lies in the admissible range. With this choice of \(\alpha\), the restricted norm on \(r\)-sets is equal to
\begin{equation}\label{eq:rNorm}
  \frac{2r(n-r)}{n}+1.
\end{equation}
Via \eqref{eq:Burkholder}, this gives the first estimate in Theorem~\ref{thm:main}.

The norm  of the full matrix \(A\) is found by taking the maximum of the block norms \eqref{eq:rNorm}. The maximum of the expression is obviously attained at \(r=n/2\). If \(n\) is even, this means that the blocks of highest norm are those related to \(n/2\)-sets, whereas for \(n\) odd, the worst bound comes from the \((n-1)/2\)- and \((n+1)/2\)-sets. The result is
\begin{equation*}
    (n/2+1)\quad\text{if }n\text{ is even},\qquad
    (n/2+1-1/2n)\quad\text{if }n\text{ is odd},
\end{equation*}
which gives the second estimate in Theorem~\ref{thm:main} via \eqref{eq:Burkholder}.

\section{Asymptotics}

Let the unit-sphere \(\Sb^{N-1}\) be equipped with its normalized rotation-invariant measure, which is denoted by \(\ud\sigma\). For an arbitrary measure space \(\Omega\), \(p\in[1,\infty]\), and \(g\in L^p(\Omega;\R^N)\), there holds
\begin{equation}\label{eq:asymp}
  \Norm{\sigma_1}{L^p(\Sb^{N-1})}\Norm{g}{L^p(\Omega;\R^N)}
  \leq\sup_{\sigma\in \Sb^{N-1}}\Norm{\sigma\cdot g}{L^p(\Omega)},
\end{equation}
which is an \(N\)-dimensional analogue of Dragi\v{c}evi\'c \& Volbergs's \cite{DraVo}, Lemma~4.1.

Indeed, by the rotation invariance,
\begin{equation*}\begin{split}
  (LHS)^p
  &=\int_{\Omega}\int_{\Sb^{N-1}}\abs{\sigma_1}^p\ud\sigma\Norm{g(x)}{\R^N}^p\ud x
   =\int_{\Omega}\int_{\Sb^{N-1}}\abs{\sigma\cdot g(x)}^p\ud x \\
  &=\int_{\Sb^{N-1}}\int_{\Omega}\abs{\sigma\cdot g(x)}^p\ud x\ud\sigma
  \leq (RHS)^p.
\end{split}\end{equation*}
The most important thing about the factor \(\Norm{\sigma_1}{L^p(\Sb^{N-1})}\) is the limiting behavior
\begin{equation}\label{eq:sigma1lim}
  \Norm{\sigma_1}{L^p(\Sb^{N-1})}\to\Norm{\sigma_1}{L^{\infty}(\Sb^{N-1})}=1\qquad\text{as }p\to\infty.
\end{equation}
Thus \eqref{eq:asymp} shows that the \(L^p\) norm of an \(\R^N\)-valued function is almost achieved by its pointwise projections to one-dimensional subspaces when \(p\) is large. This provides a tool for estimating the asymptotics of the Beurling--Ahlfors operator, indeed
\begin{equation}\label{eq:BAasymp}
  \Norm{\sigma_1}{L^p(\Sb_{\Lambda})}\Norm{S}{\bddlin(L^p(\R^n;\Lambda))}\leq
  \sup_{\sigma\in \Sb_{\Lambda}}\Norm{\sigma\cdot S}{\bddlin(L^p(\R^n;\Lambda),L^p(\R^n))},
\end{equation}
where \(\Sb_{\Lambda}\eqsim\Sb^{2^n-1}\) is the unit-sphere of the exterior algebra \(\Lambda\).

A matrix representation for \(\sigma\cdot S\) is obtained from that of \(S\); in fact, the (partial) inner product \(\sigma\cdot A\in\bddlin(\R^{2^n n},\R^n)\) will do, given a matrix \(A\) representing \(S\). Here it seems to be most convenient to use the symmetric version with all \(\alpha=\frac{1}{2}\). Using \eqref{eq:Aij}, there follows
\begin{equation*}\begin{split}
  (\sigma\cdot A)_{i,Jj}=\sum_I\sigma_I A_{Ii,Jj}
  &=\sigma_{J}\delta_{ij}(-1)^{[i\notin J]}
   +[i\in J,j\notin J]\sigma_{J\setminus i\cup j}(-1)^{\# J(i,j)} \\
  &\qquad+[i\notin J,j\in J]\sigma_{J\setminus j\cup i}(-1)^{\# J(i,j)}.
\end{split}\end{equation*}

An upper bound for the norm of such a matrix is given by the \(\ell^2\) norm of the sequence of norms of the \(n\times n\)-matrices \((\sigma\cdot A)_{i,Jj}\) for a fixed \(J\). Each of these submatrices has the form (after permutation of the indices so that all the ones in \(J\) precede those in \(J^c\))
\begin{equation*}
  (\sigma\cdot A)_J:=
  \begin{pmatrix}
     \sigma_J I_r & \Sigma_J \\
     \Sigma_J^T & -\sigma_J I_{n-r} \\
  \end{pmatrix},
\end{equation*}
where \(r=\abs{J}\) and \(\Sigma_J\) is an \(r\times(n-r)\)-matrix with entries \(\sigma_{J\setminus i\cup j}(-1)^{\# J(i,j)}\). 

Then
\begin{equation*}\begin{split}
  &\Norm{(\sigma\cdot A)_J}{\bddlin(\R^n)}^2
  =\Norm{(\sigma\cdot A)_J^T(\sigma\cdot A)_J}{\bddlin(\R^n)} \\
  &=\BNorm{
    \begin{pmatrix}
       \sigma_J^2 I_r+\Sigma_J\Sigma_J^T & 0 \\
       0 & \sigma_J^2 I_{n-r}+\Sigma_J^T\Sigma_J \\
    \end{pmatrix}}{\bddlin(\R^n)} \\
  &=\max\{\sigma_J^2+\Norm{\Sigma_J\Sigma_J^T}{\bddlin(\R^r)},
     \sigma_J^2+\Norm{\Sigma_J^T\Sigma_J}{\bddlin(\R^{n-r})}\} \\
  &=\sigma_J^2+\Norm{\Sigma_J}{\bddlin(\R^{n-r},\R^r)}^2
   \leq\sigma_J^2+\sum_{\ontop{i\in J}{j\notin J}}\sigma_{J\setminus i\cup j}^2,
\end{split}\end{equation*}
and finally
\begin{equation*}\begin{split}
  \Norm{\sigma\cdot A}{\bddlin(\R^{n 2^n},\R^n)}^2
  &\leq\sum_J\Norm{(\sigma\cdot A)_J}{\bddlin(\R^n)}^2
  \leq\sum_J\sigma_J^2+\sum_J\sum_{\ontop{i\in J}{j\notin J}}\sigma_{J\setminus i\cup j}^2 \\
  &=\sum_I\sigma_I^2\Big(1+\sum_{\ontop{i\notin I}{j\in I}}1\Big)
  =\sum_I\sigma_I^2(1+\#I\cdot\# I^c).
\end{split}\end{equation*}

The product \(\# I\cdot \# I^c\) is at most \((n/2)^2\) if \(n\) is even and \((n-1)/2\cdot(n+1)/2=(n/2)^2-1/4\) if \(n\) is odd. This proves a norm bound for the matrices \(\sigma\cdot A\) which in combination with \eqref{eq:Burkholder} yields Proposition~\ref{prop:asymp}.

\section{On the spectral multipliers of the Laplacian}

It is well known that many other operators besides the Beurling--Ahlfors transform admit the representation \eqref{eq:TvsA} and may then be estimated with the help of \eqref{eq:Burkholder}. The reader should have a look at the recent paper of Geiss, Montgomery-Smith \& Saksman~\cite{GMS}, where quite general even homogeneous Fourier multipliers (of which the Beurling--Ahlfors operator is a special case) are considered, and the results include upper as well as lower bounds for the related norms.

Here I comment briefly on the case of \eqref{eq:TvsA} when \(A=A(t)\) is a scalar-valued function (times the identity matrix, if the reader so wishes), depending only on the vertical variable. It is more natural in this context to work with complex scalars, and it is worth pointing out in any case that \eqref{eq:Burkholder} remains true in this setting, with the same proof. Then, using the Fourier transform (with the normalization \(\hat{f}(\xi)=\int f(x)e^{-i2\pi\xi\cdot x}\ud x\)),
\begin{equation*}\begin{split}
  &\int_0^{\infty}\int_{\R^n}\pair{A(t)\nabla u(x,t)}{\nabla v(x,t)}\ud x\ud t \\
  &=\int_0^{\infty}A(t)\int_{\R^n}\pair{i2\pi\xi e^{-2\pi^2\abs{\xi}^2 t}\hat{f}(\xi)}{
      i2\pi\xi e^{-2\pi^2\abs{\xi}^2 t} \hat{g}(\xi)}\ud\xi\ud t \\
  &=\int_{\R^n}\int_0^{\infty}4\pi^2\abs{\xi}^2 e^{-4\pi^2\abs{\xi}^2 t}A(t)\ud t
       \pair{\hat{f}(\xi)}{\hat{g}(\xi)}\ud\xi \\
  &=\int_{\R^n}\pair{a(4\pi^2\abs{\xi}^2)\hat{f}(\xi)}{\hat{g}(\xi)}\ud\xi
   =\int_{\R^n}\pair{a(-\triangle)f(x)}{g(x)}\ud x,
\end{split}\end{equation*}
where
\begin{equation}\label{eq:Laplace}
  a(\lambda):=\int_0^{\infty}\lambda A(t)e^{-\lambda t}\ud t.
\end{equation}

Thus the operators \(T\) corresponding to such \(A\) are spectral multipliers of the Laplacian or, equivalently, radial Fourier multipliers of the \emph{Laplace transform type}, which have been studied extensively in the literature.
The estimate \eqref{eq:Burkholder} gives in this case the norm bound
\begin{equation}\label{eq:spectral}
  \Norm{a(-\triangle)}{\bddlin(L^p(\R^n;\C^N))}\leq (p^*-1)\sup_{t\in(0,\infty)}\abs{A(t)}
\end{equation}
when \(a\) and \(A\) are related by \eqref{eq:Laplace}.

An interesting particular case consists of the imaginary powers \((-\triangle)^{is}\) which arise from \(A(t)=\Gamma(1-is)^{-1}t^{-is}\), as one readily checks from the integral representation of Euler's \(\Gamma\) function. Hence
\begin{equation}\label{eq:powers}
  \Norm{(-\triangle)^{is}}{\bddlin(L^p(\R^n;\C^N))}
  \leq\frac{(p^*-1)}{\abs{\Gamma(1-is)}}.
\end{equation}

An interesting version of the \((p^*-1)\)-principle appears in the limit as \(s\to 0\):
\begin{equation}\label{eq:limit}
  \lim_{s\to 0}\Norm{(-\triangle)^{is}}{\bddlin(L^p(\R^n;\C^N))}=p^*-1,
\end{equation}
where the existence of the limit is part of the statement. Indeed, the upper bound \(p^*-1\) for the corresponding \(\limsup\) is immediate from \eqref{eq:powers}. That the \(\liminf\) has the same lower bound is implicitly contained in the paper of Guerre-Delabri\`ere~\cite{Guerre}, as pointed out in~\cite{Aspects}. In the last-mentioned paper, the upper bound \(2(p^*-1)\) for the \(\limsup\) was obtained by using the harmonic extension method.

I conclude by mentioning that all the results \eqref{eq:spectral}, \eqref{eq:powers}, and \eqref{eq:limit} of this section extend to the generality where \(\C^N\) is replaced by any Banach space \(X\) with the \emph{unconditionality property of martingale differences} (UMD), and \(p^*-1\) by the complex UMD constant \(\beta_{p,X}^{\C}\) of the space \(X\). The number \(\beta_{p,X}^{\C}\) is defined as the smallest admissible constant in the estimate
\begin{equation*}
  \BNorm{\sum_{k=1}^N\zeta_k d_k}{L^p(\Omega;X)}\leq\beta_{p,X}^{\C}\BNorm{\sum_{k=1}^N d_k}{L^p(\Omega;X)},
\end{equation*}
which is to hold for all martingale difference sequences \((d_k)_{k=1}^N\) in \(L^p(\Omega;X)\) (the probability space \(\Omega\) and the length \(N\) also being arbitrary), and for all complex numbers \(\zeta_k\) on the unit circle. So in particular there holds
\begin{equation*}
  \lim_{s\to\infty}\Norm{(-\triangle)^{is}}{\bddlin(L^p(\R^n;X))}
  =\beta_{p,X}^{\C}.
\end{equation*}

The reason for the validity of this extension is quite simply the fact that the underlying martingale estimates are precisely the defining property of UMD spaces. I will not elaborate on this here, since the situation is sufficiently similar to the case of the harmonic extension considered in~\cite{Aspects}, and the interested reader should consult that paper for further information.

\subsection*{Acknowledgment}
Brett Wick got me interested in fine-tuning the dimensional dependence and discussed the problem with me during the thematic program ``New Trends in Harmonic Analysis'' at the Fields Institute in Toronto. The research was carried out during my stay at the Institute. I was financed by the Academy of Finland throught the project 114374 ``Vector-valued singular integrals'' and through the Finnish Centre of Excellence in Analysis and Dynamics Research.

\def\cprime{$'$}


\begin{thebibliography}{10}

\bibitem{BaMen}
R.~Ba{\~n}uelos and P.~J. M{\'e}ndez-Hern{\'a}ndez.
\newblock Space-time {B}rownian motion and the {B}eurling-{A}hlfors transform.
\newblock {\em Indiana Univ. Math. J.}, 52(4):981--990, 2003.

\bibitem{BaJan}
Rodrigo Ba{\~n}uelos and Prabhu Janakiraman.
\newblock {$L\sp p$}-bounds for the {B}eurling-{A}hlfors transform.
\newblock {\em Trans. Amer. Math. Soc.}, 360(7):3603--3612, 2008.

\bibitem{BaLin}
Rodrigo Ba{\~n}uelos and Arthur Lindeman, II.
\newblock A martingale study of the {B}eurling-{A}hlfors transform in {${\bf
  R}\sp n$}.
\newblock {\em J. Funct. Anal.}, 145(1):224--265, 1997.

\bibitem{BaWan}
Rodrigo Ba{\~n}uelos and Gang Wang.
\newblock Sharp inequalities for martingales with applications to the
  {B}eurling-{A}hlfors and {R}iesz transforms.
\newblock {\em Duke Math. J.}, 80(3):575--600, 1995.

\bibitem{Burkholder}
D.~L. Burkholder.
\newblock Boundary value problems and sharp inequalities for martingale
  transforms.
\newblock {\em Ann. Probab.}, 12(3):647--702, 1984.

\bibitem{DraVo}
Oliver Dragi{\v{c}}evi{\'c} and Alexander Volberg.
\newblock Bellman function, {L}ittlewood-{P}aley estimates and asymptotics for
  the {A}hlfors-{B}eurling operator in {$L\sp p(\Bbb C)$}.
\newblock {\em Indiana Univ. Math. J.}, 54(4):971--995, 2005.

\bibitem{GMS}
S.~Geiss, S.~Montgomery-Smith, and E.~Saksman.
\newblock On singular integral and martingale transforms.
\newblock Preprint, arXiv:math/0701516, 2007.

\bibitem{Guerre}
Sylvie Guerre-Delabri{\`e}re.
\newblock Some remarks on complex powers of {$(-\Delta)$} and {UMD} spaces.
\newblock {\em Illinois J. Math.}, 35(3):401--407, 1991.

\bibitem{Aspects}
Tuomas Hyt{\"o}nen.
\newblock Aspects of probabilistic {L}ittlewood-{P}aley theory in {B}anach
  spaces.
\newblock In B.~Randrianantoanina and N.~Randrianantoanina, editors, {\em
  Banach Spaces and their Applications in Analysis}, de Gruyter Proceedings in
  Mathematics, pages 343--355. Walter de Gruyter, Berlin, 2007.

\bibitem{IwaM}
Tadeusz Iwaniec and Gaven Martin.
\newblock Quasiregular mappings in even dimensions.
\newblock {\em Acta Math.}, 170(1):29--81, 1993.

\bibitem{Lehto}
Olli Lehto.
\newblock Remarks on the integrability of the derivatives of quasiconformal
  mappings.
\newblock {\em Ann. Acad. Sci. Fenn. Ser. A I No.}, 371:8, 1965.

\bibitem{PSW}
Stefanie Petermichl, Leonid Slavin, and Brett~W. Wick.
\newblock Heating the generalized {B}eurling-{A}hlfors operator and estimates
  for its norm.
\newblock Manuscript, 2008.

\bibitem{NaVo}
A.~Vol{\cprime}berg and F.~Nazarov.
\newblock Heat extension of the {B}eurling operator and estimates for its norm.
\newblock {\em Algebra i Analiz}, 15(4):142--158, 2003.

\end{thebibliography}

\end{document}